\documentclass[11pt]{article}
\usepackage{amsmath,amsthm}
\usepackage{amssymb,mathrsfs}
\usepackage{bm,mathtools}

\usepackage{xcolor}
\usepackage{geometry}
\usepackage[colorlinks=true,
linkcolor=blue,citecolor=blue,
urlcolor=blue]{hyperref}

\usepackage{microtype}

\theoremstyle{plain}
\newtheorem{theorem}{Theorem}

\theoremstyle{definition}
\newtheorem{remark}{Remark}

\allowdisplaybreaks
\DeclareMathOperator{\trans}{\mathrm{T}}

\title{A short proof of an inequality on spectral radius perturbation}

\author{Lele Liu\footnote{School of Mathematical Sciences, Anhui University, Hefei 230601, 
P.R. China. E-mail: \texttt{liu@ahu.edu.cn} (L. Liu). Supported by the National
Nature Science Foundation of China (No. 12471320) and the Anhui Provincial Natural Science Foundation (No. 2408085Y003).}
~~ and ~~
Bo Ning\footnote{College of Computer Science, Nankai University, Tianjin 300350, P.R. China.
E-mail: \texttt{bo.ning@nankai.edu.cn} (B. Ning). Partially supported by the National Nature Science
Foundation of China (No. 12371350).}}

\date{}

\begin{document}
\maketitle

\begin{abstract}
Let $G$ be a simple graph, and denote by $\lambda(G)$ its spectral radius.  
Sun and Das (2020) established that for any non-isolated vertex $v$ with degree $d(v)$,
\[
\lambda(G)\leq \sqrt{\lambda(G-v)^2 + 2d(v) - 1},
\]
which is a conjecture original posed by Guo, Wang, and Li (2019). 
Sun and Das's proof use several tools from spectral graph theory. In this short note, 
we provide a concise and self-contained proof of this inequality using matrix analysis.
\par\vspace{2mm}

\par\vspace{2mm}

\end{abstract}

\section{Introduction}

Let $G$ be a simple graph and let $\lambda(G)$ denote the spectral radius of 
the adjacency matrix of $G$. It is well-known that $\lambda(G-v)<\lambda(G)$ for any connected graph $G$ and any vertex $v\in V(G)$ (see \cite[Proposition~1.3.9]{Cvetkovic-Rowlinson-Simic2010}).
A natural problem in spectral graph theory is to ask for the estimate of $\lambda(G)-\lambda(G-v)$ for any graph $G$ and a vertex $v\in V(G)$. One can see that the difference can tend to infinite and even tend to 0. Consider $G=K_{1,n-1}$ and $v$ is the center or any other leaf-vertex.
One classical inequality is that $\lambda(G)\leq \lambda(G-v)+1$ if $v$ is a minimum degree vertex.
Nikiforov \cite[Lemma~11]{Nikiforov2010} proved that $\lambda(G)\leq \lambda(G-v_k)\frac{1-x^2_k}{1-2x^2_k}$, where $x_k$
is the smallest eigen-component of a unit eigenvector and $v_k$ is the vertex with respect to $x_k$.

Motivated by Hong's famous inequality \cite{Hong1993}, Guo, Wang, and Li \cite{Guo-Wang-Li2019}  conjectured that for any simple graph $G$ and
non-isolated vertex $v$,
\[
\lambda(G)\leq \sqrt{\lambda(G-v)^2+2d(v)-1},
\]
and proved it when $d(v)=1$. 
Guo et al. \cite{Guo-Wang-Li2019} observed that if it were true, then it would imply Hong's inequality. Recall that Hong's inequality states that $\lambda(G)\leq\sqrt{2m-n+1}$ if $\delta(G)\geq 1$, where $m$ is the number of edges, and $n$ is the number of vertices of $G$ respectively. Later, 
Sun and Das \cite{Sun-Das2020} 
confirmed the conjecture. 
This spectral inequality has found several useful applications in spectral graph theory, see, for example, advancement on the spectral consecutive cycle problem \cite{Li-Ning2023}, a complete characterization of spectral extremal graphs on $P_k$-free graphs \cite{Ai-Lei-Ning-Shi2026+}, and \cite{Zhai-Lin-Shu2021}. For a powerful extension, see \cite{Jin-Zhang-Zhang2024}.

\begin{theorem}[\cite{Sun-Das2020}]\label{Thm:GuoetalConj}
Let $G$ be a simple graph and let $v\in V(G)$ be a non-isolated vertex of degree $d(v)$. Then
$$
\lambda(G)\leq \sqrt{\lambda(G-v)^2 + 2d(v) - 1}.
$$
If $G$ is connected, the equality holds if and only if $G\cong K_{1,n-1}$ and 
$d(v) = 1$ or $G\cong K_n$. 
\end{theorem}

The purpose of this note is to provide a short and self-contained proof of Theorem \ref{Thm:GuoetalConj},
with aid of matrix theory. 
There are several beautiful works which use matrix theory to attack open problems on graph theory, see, for example, Ando and Lin's paper \cite{Ando-Lin2015}, solving a conjecture of Wocjan and Elphick \cite{Wocjan-Elphick2013}.

Let $A$ and $B$ be two symmetric matrices. We write $A \succeq B$ if $A-B$ is positive 
semi-definite. For graph notation and terminology undefined here, we refer the reader to \cite{Bondy-Murty2008}.

\section{A short proof of Theorem \ref{Thm:GuoetalConj}}
Set $\lambda:=\lambda(G)$ and $\mu:=\lambda(G-v)$ for brevity. Clearly, $\lambda\geq\mu$. 
If $\lambda = \mu$, the result follows immediately. Thus, we may assume $\lambda > \mu$.

Relabeling the vertices of $G$ if necessary, we may assume that the adjacency matrix 
$A(G)$ has the block form
\[
A(G)= 
\begin{bmatrix}
0 & \bm{b}^{\mathrm{T}} \\
\bm{b} & B
\end{bmatrix},
\]
where $B=A(G-v)$ and $\bm{b}$ is the $0$--$1$ column vector indicating the neighbors of $v$ 
in $G-v$. In particular, $\bm{b}^{\trans} \bm{b} = d(v)$.

Since $\lambda>\mu=\lambda(B)$, the matrix $\lambda I-B$ is positive definite and hence invertible. 
Moreover, $\det\,(\lambda I-A(G))=0$. Taking the Schur complement of $\lambda I - B$ in $\lambda I - A(G)$, 
we obtain
\[
0 = \det\,(\lambda I - B)\cdot \big(\lambda - \bm{b}^{\trans} (\lambda I - B)^{-1} \bm{b} \big).
\]
It follows from $\det\,(\lambda I - B) \neq 0$ that
$\lambda = \bm{b}^{\trans} (\lambda I - B)^{-1} \bm{b}$.
Now, let $t$ be any eigenvalue of $B$. Since $|t|\leq \mu < \lambda$, we have
\[
\frac{\lambda+t}{\lambda^2-\mu^2}-\frac{1}{\lambda-t} =
\frac{\mu^2-t^2}{(\lambda^2-\mu^2)(\lambda-t)}\geq 0
\implies \frac{\lambda + t}{\lambda^2 - \mu^2} \geq \frac{1}{\lambda - t}.
\]
We immediately obtain $\frac{\lambda I+B}{\lambda^2-\mu^2} \succeq (\lambda I-B)^{-1}$.
Indeed, let $B = U\,\mathrm{diag}(t_1,t_2,\ldots,t_{n-1})\,U^{\trans}$, 
where $U$ is orthogonal. Then
\[
M: = \frac{\lambda I {+} B}{\lambda^2 {-} \mu^2} - (\lambda I {-} B)^{-1}
= U\,\mathrm{diag}\!\left(\frac{\lambda {+} t_1}{\lambda^2 {-} \mu^2} {-} 
\frac{1}{\lambda {-} t_1},\ldots,\frac{\lambda {+} t_{n-1}}{\lambda^2 {-} \mu^2} 
{-} \frac{1}{\lambda {-} t_{n-1}}\right) U^{\trans}
\succeq 0.
\]
Recall that $M\succeq 0$ and $\lambda = \bm{b}^{\trans} (\lambda I - B)^{-1} \bm{b}$. Thus, we obtain
\begin{equation}\label{eq:firstbound}
\frac{\lambda\, d(v) + \bm{b}^{\trans} B \bm{b}}{\lambda^2 - \mu^2}
- \lambda = \frac{\lambda\, \bm{b}^{\trans} \bm{b} + \bm{b}^{\trans} B \bm{b}}{\lambda^2 - \mu^2} 
- \bm{b}^{\trans} (\lambda I - B)^{-1} \bm{b} = \bm{b}^{\trans} M\bm{b} \geq 0.
\end{equation}

To finish the proof, we shall estimate $\bm{b}^{\trans} B \bm{b}$. 
Observe that $\bm{b}^{\trans} B \bm{b} = 2e\big(G[N(v)]\big)$, 
where $G[N(v)]$ is the subgraph induced by the neighborhood of $v$.
Set $m:= \bm{b}^{\trans} B \bm{b}$ for short. If $d(v) = 1$, then clearly $m=0$, 
and hence $\lambda^2 - \mu^2 \leq d(v) = 1 = 2d(v) - 1$ by \eqref{eq:firstbound}. 
Next, we may therefore assume that $d(v)\geq 2$. Let $F=G[\{v\}\cup N(v)]$ be the 
subgraph induced by the closed neighborhood of $v$. Then $|V(F)| = d(v) + 1$ 
and $2e(F)= 2 d(v) + m$. Thus, the average degree of $F$ is
$\overline d(F) = \frac{2d(v) + m}{d(v) + 1}$.
Since the spectral radius of a graph is at least its average degree, we have 
$\lambda\geq \lambda(F)\geq \overline d(F)$. Moreover, 
since $m = 2e(G[N(v)])\leq d(v) (d(v) - 1)$, we obtain
\begin{equation}\label{eq:average-deg-F}
\frac{2d(v) + m}{d(v) + 1} - \frac{m}{d(v) - 1}
= \frac{2\big(d(v) (d(v) - 1) - m\big)}{d(v)^2 - 1} \geq 0.
\end{equation}
It follows that $\lambda\geq \overline d(F)\geq \frac{m}{d(v) - 1}$, and hence
$m = \bm{b}^{\trans} B \bm{b} \leq \lambda\, (d(v) - 1)$.
Substituting this into \eqref{eq:firstbound} yields
\[
\lambda\leq \frac{\lambda\, d(v) + \lambda\, (d(v) - 1)}{\lambda^2 - \mu^2} 
= \frac{\lambda\, (2d(v) - 1)}{\lambda^2 - \mu^2}.
\]
This rearranges to $\lambda^2 - \mu^2\leq 2d(v) - 1$, which is exactly the desired inequality.

Finally, we consider the equality case. The sufficiency is clear,  so we focus on proving the necessity. 
If $d(v) = 1$, let $u$ be the unique neighbor of $v$. Then $\bm{b}=e_u$, where $e_u$ is the indicator vector of $u$. By \eqref{eq:firstbound}, 
we have $e_u^{\trans} M e_u = 0$, which implies $M e_u = \bm{0}$. Consequently, $(\lambda I - B) Me_u = \bm{0}$.
On the other hand, 
\[
(\lambda I - B) M = \frac{(\lambda I - B) (\lambda I + B)}{\lambda^2 - \mu^2} - I 
= \frac{\mu^2 I - B^2}{\lambda^2 - \mu^2}.
\]
It follows that $B^2e_u = \mu^2 e_u$. Observe that for each $v\in V(H)$, 
$(B^2e_u)_v = |N_H(v)\cap N_H(u)|$. Thus, for each $v\neq u$,
$|N_H(v)\cap N_H(u)| = 0$. This means that there are no two neighbors of $u$ adjacent, 
and there is no vertex at distance $2$ from $u$ in $H$. Since $H$ is connected, 
it follows that $H$ is a star with the center $u$. If $d(v)\geq 2$, then necessarily 
$\lambda = \lambda(F)$, and the equality must hold in \eqref{eq:average-deg-F}. 
Since $\lambda = \lambda (F)$ and $G$ is connected, we must have $F\cong G$. 
Moreover, equality in \eqref{eq:average-deg-F} implies that $m=2e(G[N(v)])=d(v) (d(v)-1)$.
Thus, $G$ is a complete graph.
The proof is complete.

\begin{remark}
We point out the following interesting corollary of Theorem 1:
Choose $v$ such that $d(v)\leq \frac{2m}{n}$. Then 
$$
\lambda(G)\leq \lambda(G-v)+1.
$$
Indeed, $d(v)\leq \frac{2m}{n}\leq \lambda(G)$. The inequality follows.
\end{remark}

\end{document}